\documentclass[mlq,a4paper]{w-art}
\usepackage{times,cite,w-thm}

\usepackage{amssymb}

\newcommand{\fa}{\forall}

\newcommand{\Si}{\Sigma}

\newcommand{\Sio}{\Sigma^\omega}
\newcommand{\ra}{\longrightarrow}
\newcommand{\hs}{\hspace{12mm}

\noi}
\newcommand{\lra}{\leftrightarrow}

\newcommand{\ite}{\item}

\newcommand{\ol}{ $\omega$-language}

\newcommand{\om}{\omega}
\newcommand{\nl}{\newline}
\newcommand{\noi}{\noindent}
\newcommand{\tla}{\twoheadleftarrow}

\newcommand{\bormxi}{{\bf\Pi}^{0}_{\xi}}

\newcommand{\boraone}{{\bf\Sigma}^{0}_{1}}

\newcommand{\boraxi}{{\bf\Sigma}^{0}_{\xi}}

\begin{document}

\DOIsuffix{theDOIsuffix}

\Volume{}
\Month{}
\Year{2009}

\pagespan{1}{10}

\Receiveddate{XXXX}
\Reviseddate{XXXX}
\Accepteddate{XXXX}
\Dateposted{XXXX}



\title[On  some sets of dictionaries whose $\omega$-powers  have a given complexity]{On  some sets 
of dictionaries whose $\omega$-powers \\ have a given complexity}

\author[O. Finkel]{Olivier Finkel
  \footnote{Corresponding author\quad E-mail:~\textsf{finkel@logique.jussieu.fr}}}
\address{Equipe de Logique Math\'ematique\\ CNRS et  Universit\'e Paris Diderot Paris 7
 \\ UFR de Math\'ematiques case 7012, site Chevaleret,\\75205 Paris Cedex 13,  France}

\begin{abstract}
\noi A dictionary  is a set of finite words over some finite alphabet $X$. 
The $\omega$-power of a dictionary $V$ is the set of infinite words  obtained by infinite concatenation of words in $V$. 
 Lecomte studied in \cite{Lecomte-JSL} the complexity of the set of dictionaries whose associated $\om$-powers have a given complexity. 
In particular,  he considered the sets $\mathcal{W}({\bf\Si}^0_{k})$ (respectively, $\mathcal{W}({\bf\Pi}^0_{k})$, $\mathcal{W}({\bf\Delta}_1^1)$) 
of dictionaries $V \subseteq 2^\star$ whose $\om$-powers are ${\bf\Si}^0_{k}$-sets (respectively, ${\bf\Pi}^0_{k}$-sets, Borel sets). 
In this paper we first establish a new relation between the sets  $\mathcal{W}({\bf\Sigma}^0_{2})$ and 
$\mathcal{W}({\bf\Delta}_1^1)$,  
showing that the set  $\mathcal{W}({\bf\Delta}_1^1)$ is ``more complex" than the set $\mathcal{W}({\bf\Sigma}^0_{2})$. 
As an application we improve the lower bound on the complexity of   $\mathcal{W}({\bf\Delta}_1^1)$  given by Lecomte, showing that 
$\mathcal{W}({\bf\Delta}_1^1)$ is in ${\bf\Sigma}^1_2(2^{2^{\star}})\!\setminus {\bf\Pi}^0_2$. Then we prove that, for every integer $k\geq 2$, 
(respectively, $k\geq 3$)
 the set of dictionaries
$\mathcal{W}({\bf\Pi}^0_{k+1})$  (respectively, $\mathcal{W}({\bf\Si}^0_{k+1})$)     is ``more complex" than the set of dictionaries 
$\mathcal{W}({\bf\Pi}^0_{k})$ (respectively, $\mathcal{W}({\bf\Si}^0_{k})$) . 
\end{abstract}

\maketitle 

\section{Introduction}

A finitary language, called here also a dictionary as in \cite{Lecomte-JSL},  is a set of finite words over some finite alphabet $X$. 
The $\omega$-power of a dictionary $V$ is the set of infinite words  obtained by infinite concatenation of words in $V$. 
The  $\omega$-powers  appear very naturally in Theoretical Computer Science and in Formal Language Theory,   in the characterization of the  classes 
of  languages of infinite words accepted by finite automata or  by pushdown automata, \cite{Staiger97}.  

\hs Since the set  of infinite words over a finite alphabet $X$ is usually equipped 
with the Cantor topology, the  question of the topological complexity of the $\om$-powers of 
finitary  languages naturally arises.  It has been posed by 
Niwinski \cite{Niwinski90},  Simonnet \cite{Simonnet92} and  Staiger \cite{Staiger97b}. 
\nl Firstly it is easy to see that the  $\om$-power of a finitary language $V$  is always an analytic set because 
it is  the continuous image of either  a compact set $\{1, \ldots ,n\}^\om$ for $n\geq 0$, 
or the Baire space $\om^\om$.  
\nl The first example of a finitary language $L$ such that the $\om$-power $L^\om$ is analytic but not Borel, 
and even ${\bf \Si}_1^1$-complete, was obtained in  \cite{Fin03a}. 
Amazingly the language $L$ has a very simple description 
and  was obtained   via a coding of the infinite labelled  binary trees. The construction will be recalled below. 
For the Borel $\om$-powers, after some partial results obtained in  \cite{Fin01a,Fin04-FI,Fin-Dup}, the 
question of the Borel hierarchy of  $\om$-powers of finitary languages has been solved  recently by Finkel and Lecomte in \cite{Fin-Lec}, where 
a very surprising result is proved, showing that actually $\om$-powers exhibit a great topological complexity. For every non-null countable ordinal 
$\alpha$ there exist some ${\bf \Si}^0_\alpha$-complete $\om$-powers and also some ${\bf \Pi}^0_\alpha$-complete $\om$-powers.

\hs  Another question naturally arises about $\om$-powers and descriptive set theory.  It  has been firstly studied by Lecomte in 
\cite{Lecomte-JSL}.  He asked about  the complexity of the set of dictionaries whose associated $\om$-powers have a given complexity. 
The set $\mathcal{W}({\bf\Sigma}^0_{\xi})$ (respectively,  $\mathcal{W}({\bf\Pi}^0_{\xi})$, $\mathcal{W}({\bf\Delta}_1^1)$) 
is the set of dictionaries over the alphabet $2=\{0, 1\}$ 
whose $\om$-powers are  ${\bf\Sigma}^0_{\xi}$-sets
(respectively,  ${\bf\Pi}^0_{\xi}$-sets,  Borel sets). The set of dictionaries over the alphabet $2=\{0, 1\}$ can be  naturally equipped with the Cantor topology. 
Then Lecomte proved  that $\mathcal{W}({\bf\Sigma}^0_{2})$ is in ${\bf\Sigma}^1_2(2^{2^{\star}})\!\setminus {\bf\Pi}^0_2$ and that 
all the other sets $\mathcal{W}({\bf\Sigma}^0_{\xi})$, $\mathcal{W}({\bf\Pi}^0_{\xi})$, and $\mathcal{W}({\bf\Delta}_1^1)$ 
 are in ${\bf\Sigma}^1_2(2^{2^{\star}})\!\setminus\! D_2(\boraone )$, where 
 $D_2(\boraone )$ is the class of $2$-differences of open sets, that  is,   the class of sets which are intersections of an open set and of a closed set. 
It is proved  in \cite{Fink-Lec2} that for each countable ordinal $\xi \geq 3$ the sets $\mathcal{W}({\bf\Sigma}^0_{\xi})$ and 
$\mathcal{W}({\bf\Pi}^0_{\xi})$
are actually 
${\bf\Pi}_{1}^1$-hard. 
In this paper we obtain first a new relation between the sets $\mathcal{W}({\bf\Sigma}^0_{2})$ and $\mathcal{W}({\bf\Delta}_1^1)$, 
showing that $\mathcal{W}({\bf\Sigma}^0_{2})$ is continuously reducible to $\mathcal{W}({\bf\Delta}_1^1)$, which means 
that the set  $\mathcal{W}({\bf\Delta}_1^1)$ is ``more complex" than the set $\mathcal{W}({\bf\Sigma}^0_{2})$. 
As an application we improve the lower bound on the complexity of   $\mathcal{W}({\bf\Delta}_1^1)$  given by Lecomte, showing that 
$\mathcal{W}({\bf\Delta}_1^1)$ is in ${\bf\Sigma}^1_2(2^{2^{\star}})\!\setminus {\bf\Pi}^0_2$. Then we prove that, for every integer $k\geq 2$,
(respectively, $k\geq 3$)
 the set of dictionaries
$\mathcal{W}({\bf\Pi}^0_{k+1})$  (respectively, $\mathcal{W}({\bf\Si}^0_{k+1})$)   
  is ``more complex" than the set of dictionaries $\mathcal{W}({\bf\Pi}^0_{k})$ (respectively, $\mathcal{W}({\bf\Si}^0_{k})$) .

\hs  The paper is organized as follows. In Section 2 we recall some notations of formal language theory and some notions of topology. We prove our results 
in Section 3.  Some concluding remarks are given in Section 4.

\section{Borel and projective hierarchies}

We  use usual notations of formal language theory which may be found for instance in \cite{Thomas90,Staiger97}. 
\nl  When $X$ is a finite alphabet, a {\it non-empty finite word} over X is any 
sequence $x=a_1\ldots a_k$, where $a_i\in X$ 
for $i=1,\ldots ,k$ , and  $k$ is an integer $\geq 1$. The {\it length}
 of $x$ is $k$, denoted by $|x|$.
 The {\it empty word} has no letter and is denoted by $\lambda$; its length is $0$. 
 $X^\star$  is the {\it set of finite words} (including the empty word) over $X$, and $X^+=X^\star \setminus \{\lambda\}$ is the set of 
 {\it non-empty finite words}. 
A {\it   finitary language}, called here also a {\it dictionary}, 
over the alphabet $X$   is a subset of  $X^\star$.

 \hs An $\om$-{\it word} over $X$ is an $\om$ -sequence $a_1 \ldots a_n \ldots$, where for all 
integers $ i\geq 1$, ~
$a_i \in X$.  When $\sigma$ is an $\om$-word over $X$, we write
 $\sigma =\sigma(1)\sigma(2)\ldots \sigma(n) \ldots $,  where for all $i$,~ $\sigma(i)\in X$,
and $\sigma[n]=\sigma(1)\sigma(2)\ldots \sigma(n)$  for all $n\geq 1$ and $\sigma[0]=\lambda$.
\nl 
 The usual concatenation product of two finite words $u$ and $v$ is 
denoted $u\cdot v$ (and sometimes just $uv$). This product is extended to the product of a 
finite word $u$ and an $\om$-word $v$: the infinite word $u\cdot v$ is then the $\om$-word such that:
\nl $(u\cdot v)(k)=u(k)$  if $k\leq |u|$ , and 
 $(u\cdot v)(k)=v(k-|u|)$  if $k>|u|$.
\nl   The {\it prefix relation} is denoted $\sqsubseteq$: a finite word $u$ is a {\it prefix} 
of a finite word $v$ (respectively,  an infinite word $v$), denoted $u\sqsubseteq v$,  
 if and only if there exists a finite word $w$ 
(respectively,  an infinite word $w$), such that $v=u\cdot w$.
\nl   
 The {\it set of } $\om$-{\it words} over  the alphabet $X$ is denoted by $X^\om$.
An  $\om$-{\it language} over an alphabet $X$ is a subset of  $X^\om$.  
\nl  We shall denote  $X^{\leq \om}=X^\star \cup X^\om$ the set of {\it finite or infinite } words over the alphabet $X$. 

\hs  We assume the reader to be familiar with basic notions of topology which
may be found in \cite{Moschovakis80,LescowThomas,Kechris94,Staiger97,PerrinPin}.
There is a natural metric on the set $X^\om$ of  infinite words 
over a finite alphabet 
$X$ containing at least two letters. It is  called the {\it prefix metric} and is defined as follows. For $u, v \in X^\om$ and 
$u\neq v$ let $\delta(u, v)=2^{-l_{\mathrm{pref}(u,v)}}$ where $l_{\mathrm{pref}(u,v)}$ 
 is the first integer $n$
such that the $(n+1)^{st}$ letter of $u$ is different from the $(n+1)^{st}$ letter of $v$. 
This metric induces on $X^\om$ the usual  Cantor topology for which the {\it open subsets} of 
$X^\om$ are of  the form $W\cdot X^\om$, where $W\subseteq X^\star$.
A set $L\subseteq X^\om$ is a {\it closed set} iff its complement $X^\om \setminus L$ 
is an open set.
Define now the {\it Borel Hierarchy} of subsets of $X^\om$:

\begin{defn}
For a non-null countable ordinal $\alpha$, the classes ${\bf \Si}^0_\alpha$
 and ${\bf \Pi}^0_\alpha$ of the Borel Hierarchy on the topological space $X^\om$ 
are defined as follows:
\nl ${\bf \Si}^0_1$ is the class of open subsets of $X^\om$, 
 ${\bf \Pi}^0_1$ is the class of closed subsets of $X^\om$, 
\nl and for any countable ordinal $\alpha \geq 2$: 
\nl ${\bf \Si}^0_\alpha$ is the class of countable unions of subsets of $X^\om$ in 
$\bigcup_{\gamma <\alpha}{\bf \Pi}^0_\gamma$.
 \nl ${\bf \Pi}^0_\alpha$ is the class of countable intersections of subsets of $X^\om$ in 
$\bigcup_{\gamma <\alpha}{\bf \Si}^0_\gamma$.
\end{defn}

\noi For 
a countable ordinal $\alpha$,  a subset of $\Si^\om$ is a Borel set of {\it rank} $\alpha$ iff 
it is in ${\bf \Si}^0_{\alpha}\cup {\bf \Pi}^0_{\alpha}$ but not in 
$\bigcup_{\gamma <\alpha}({\bf \Si}^0_\gamma \cup {\bf \Pi}^0_\gamma)$.

\hs    
There exists another hierarchy beyond the Borel hierarchy, which is called the 
projective hierarchy. The classes ${\bf \Si}^1_n$ and ${\bf \Pi}^1_n$, for integers $n\geq 1$, 
of the projective hierarchy are  obtained from  the Borel hierarchy by 
successive applications of operations of projection and complementation.
The first level of the projective hierarchy consists of  the class of {\it analytic sets}, and the class of {\it co-analytic sets} which are complements of 
analytic sets.  
In particular,  
the class of Borel subsets of $X^\om$ is strictly included in 
the class  ${\bf \Si}^1_1$ of {\it analytic sets}.  The class of analytic sets is also the class of the 
continuous images  of Borel sets. 

\hs We now recall the notion of Wadge reducibility, which will be fundamental in  the sequel. 

\begin{defn}[Wadge \cite{Wadge83}] Let $X$, $Y$ be two finite alphabets. 
For $L\subseteq X^\om$ and $L'\subseteq Y^\om$, $L$ is said to be Wadge reducible to $L'$
($L\leq _W L')$ iff there exists a continuous function $f: X^\om \ra Y^\om$, such that
$L=f^{-1}(L')$.
\nl $L$ and $L'$ are Wadge equivalent iff $L\leq _W L'$ and $L'\leq _W L$. 
This is  denoted by $L\equiv_W L'$. 
\end{defn}

\noi
 The relation $\leq _W $  is reflexive and transitive,
 and $\equiv_W $ is an equivalence relation.
\nl The {\it equivalence classes} of $\equiv_W $ are called {\it Wadge degrees}. 
\nl  For $L\subseteq X^\om$ and $L'\subseteq Y^\om$, if   
$L\leq _W L'$ and $f$ is a continuous 
function from $ X^\om$  into $Y^\om$ with 
$L=f^{-1}(L')$,   then $f$ is called a continuous reduction of $L$ to 
$L'$. Intuitively it means that $L$ is less complicated than $L'$ because 
to check whether $x\in L$ it suffices to check whether $f(x)\in L'$ where $f$ 
is a continuous function.

\hs  Recall that each  Borel class ${\bf \Si}^0_\alpha$ and ${\bf \Pi}^0_\alpha$  is closed under inverse images by continuous functions 
and that a  set $L\subseteq X^\om$ is a ${\bf \Si}^0_\alpha$
 (respectively ${\bf \Pi}^0_\alpha$)-{\it complete set} iff for any set 
$L'\subseteq Y^\om$, $L'$ is in 
${\bf \Si}^0_\alpha$ (respectively ${\bf \Pi}^0_\alpha$) iff $L'\leq _W L $.

\hs  There is a close relationship between Wadge reducibility
 and games that  we now introduce.  

\begin{defn} Let 
$L\subseteq X^\om$ and $L'\subseteq Y^\om$. 
The Wadge game  $W(L, L')$ is a game with perfect information between two players. 
Player 1  is in charge of $L$ and Player 2  is in charge of $L'$.
\nl Player 1 first writes a letter $a_1\in X$, then Player 2 writes a letter
$b_1\in Y$, then Player 1 writes a letter $a_2\in  X$, and so on. 
\nl The two players alternatively write letters $a_n$ of $X$ for Player 1 and $b_n$ of $Y$
for Player 2.
\nl After $\om$ steps, Player 1 has written an $\om$-word $a\in X^\om$ and Player 2
has written an $\om$-word $b\in Y^\om$.
 Player 2 is allowed to skip, even infinitely often, provided he really writes an
$\om$-word in  $\om$ steps.
\nl Player 2 wins the play iff [$a\in L \lra b\in L'$], i.e. iff : 
\begin{center}
  [($a\in L ~{\rm and} ~ b\in L'$)~ {\rm or} ~ 
($a\notin L ~{\rm and}~ b\notin L'~{\rm and} ~ b~{\rm is~infinite}  $)].
\end{center}
\end{defn}

\noi
Recall that a strategy for Player 1 is a function 
$\sigma :(Y\cup \{s\})^\star\ra X$.
And a strategy for Player 2 is a function $f:X^+\ra Y\cup\{ s\}$.
\nl A strategy $\sigma$ is a winning stategy  for Player 1 iff he always wins the play when
 he uses the strategy $\sigma$, i.e. when the  $n^{th}$  letter he writes is given
by $a_n=\sigma (b_1\cdots b_{n-1})$, where $b_i$ is the letter written by Player 2 
at the step $i$ and $b_i=s$ if Player 2 skips at the step $i$.
 A winning strategy for Player 2 is defined in a similar manner.

\hs   Martin's Theorem states that every Gale-Stewart game $G(B)$,  where  $B$ is a Borel set, 
is determined,  see \cite{Kechris94}. This implies the following determinacy result:

\begin{theorem} [Wadge] Let $L\subseteq X^\om$ and $L'\subseteq Y^\om$ be two Borel sets, where
$X$ and $Y$ are finite  alphabets. Then the Wadge game $W(L, L')$ is determined:
one of the two players has a winning strategy. And $L\leq_W L'$ iff Player 2 has a 
winning strategy  in the game $W(L, L')$.
\end{theorem}

\section{$\om$-powers and sets of dictionaries}

\noi  Recall that, for $V\subseteq X^\star$, the $\om$-language 
\[V^\om = \{ u_1\cdot u_2  \cdots u_n \cdots 
 \mid \fa i\geq 1 ~ u_i \in V\setminus \{\lambda\}  \}\]
\noi  is the $\om$-power of the language, or dictionary, $V$.

\hs A dictionary over the alphabet $X$ may be seen as an element of the space $2^{X^\star}$, i.e. the set of functions from  $X^\star$ into $2$, 
where $2=\{0,1\}$ is a two letter alphabet. 
 The space $2^{X^\star}$ is naturally equipped with 
the product topology of the discrete topology on $2=\{0,1\}$. 
The set $X^\star$ 
of finite words over the alphabet $X$ is countable so there is a bijection between $X^\star$ and $\om$ and the topological space $2^{X^\star}$ is in fact 
homeomorphic to the Cantor space $2^\om$. The notions of Borel and projective hierarchies on the space $2^{X^\star}$ are obtained in the 
same way as above in the case of the Cantor space $X^\om$.  

\hs Lecomte introduced in \cite{Lecomte-JSL} the following sets of dictionaries.  For a non null countable ordinal $\xi$, we set  

\[\begin{array}{ll}
\mathcal{W}({\bf\Sigma}^0_{\xi})\!\!\!\! & :=\!\{ A\!\subseteq\! 2^{\star}\mid A^{\om} \mbox{ is a } \boraxi \mbox{-set}~ \}
\hbox{\rm ,}\cr & \cr
\mathcal{W}({\bf\Pi}^0_{\xi})\!\!\!\! & :=\!\{ A\!\subseteq\! 2^{\star}\mid A^{\om} \mbox{ is a }  \bormxi \mbox{-set} ~\}\hbox{\rm ,}\cr & \cr
\ \mathcal{W}({\bf\Delta}_1^1) \!\!\!\! & :=\!\{ A\!\subseteq\! 2^{\star}\mid A^{\om} \mbox{ is a Borel set } \}.
\end{array}\]

\noi Lecomte  proved in  \cite{Lecomte-JSL} that $\mathcal{W}({\bf\Sigma}^0_{2})$ is in 
${\bf\Sigma}^1_2(2^{2^{\star}})\!\setminus {\bf\Pi}^0_2$ and that 
all the other sets $\mathcal{W}({\bf\Sigma}^0_{\xi})$, $\mathcal{W}({\bf\Pi}^0_{\xi})$, and $\mathcal{W}({\bf\Delta}_1^1)$ are 
in ${\bf\Sigma}^1_2(2^{2^{\star}})\!\setminus\! D_2(\boraone )$, where  
 $D_2(\boraone )$ is the class of $2$-differences of open sets, that is,  the class of sets which are intersections of an open set and of a closed set. 
Finkel and Lecomte showed in \cite{Fink-Lec2} that for each countable ordinal $\xi \geq 3$ the sets $\mathcal{W}({\bf\Sigma}^0_{\xi})$ 
and $\mathcal{W}({\bf\Pi}^0_{\xi})$ are actually 
${\bf\Pi}_{1}^1$-hard. This gives a much better lower bound on  the complexity of these sets,  but their complexity is not completely determined. 

\hs Staiger gave in  \cite{Staiger97b}  a  characterization of   the set  $\mathcal{W}({\bf\Sigma}^0_{1})$  (respectively, $\mathcal{W}({\bf\Pi}^0_{1})$). 
He gave  in \cite{Staiger97b} an example of a dictionary $V \in \mathcal{W}({\bf\Sigma}^0_{1} \setminus {\bf\Pi}^0_{1})$, and 
also an example of a $W \in \mathcal{W}({\bf\Delta}^0_{2}) \setminus \mathcal{W}({\bf\Sigma}^0_{1}  \cup {\bf\Pi}^0_{1})$. We refer the reader to 
\cite{Fin-Lec, Fink-Lec2} for an example of a $W \in \mathcal{W}({\bf\Sigma}^0_{2} \setminus  {\bf\Pi}^0_{2})$. 

\hs   In this paper we show that the set  $\mathcal{W}({\bf\Delta}_1^1)$ is more complex than the set  $\mathcal{W}({\bf\Sigma}^0_{2})$. 
As an application we improve the lower bound on the  complexity of the 
set $\mathcal{W}({\bf\Delta}_1^1)$. 

\hs We have already mentioned in the introduction the existence of  a dictionary $L$ such that $L^\om$ is ${\bf \Si}_1^1$-complete, and  hence non Borel. 
We now give a simple  construction of such a  language $L$ 
using the notion of substitution that we now recall, (see \cite{Fin03a} for more details). 

\hs  A {\it substitution}  is defined by a mapping 
$f: X\ra \mathcal{P}(Y^\star)$, where $X =\{a_1, \ldots ,a_n\}$  and $Y$ are two finite alphabets.  
For each integer $i\in [1;n]$, $f(a_i)=L_i$ is a finitary language 
over the alphabet  $Y$.
\nl Now this mapping is extended in the usual manner to finite words:
$f(a_{i_1} \cdots a_{i_n})= L_{i_1} \cdots L_{i_n}$,   
and to finitary languages $L\subseteq X^\star$: 
$f(L)=\cup_{x\in L} f(x)$.  
\nl If for each integer $i\in [1;n]$ the language  $L_i$ does not 
contain the empty word, then the mapping $f$ may be extended to $\om$-words: ~~
\[f(x(1)\cdots x(n)\cdots )= \{u_1\cdots u_n \cdots  \mid \fa i \geq 1 \quad u_i\in f(x(i))\}\] 
\noi and to \ol s $L\subseteq X^\om$ by setting  $f(L)=\cup_{x\in L} f(x)$.   

\hs 
Now let $X=\{0,1\}$,    $d$ be a new letter not in $X$, and 
\[D=\{ u\cdot d\cdot v  \mid  u, v \in X^\star ~and~ [ ( |v|=2|u|)~~ or ~~( |v|=2|u|+1) ]~ 
\}\] 

\noi Let $g: X \ra \mathcal{P}((X \cup \{d\})^\star)$ be the substitution 
defined by 
$g(a)=a\cdot D$. 

\hs  Notice that if $V^\om$ is an $\om$-power then   
$g(V^\om)=(g(V))^\om$ is also an $\om$-power.

\hs If 
$W=0^\star\cdot 1$ then 
 $W^\om=(0^\star\cdot 1)^\om$ is the set of $\om$-words over the alphabet $X$ 
containing  infinitely many occcurrences of the letter $1$.
It is a well known  example of an $\om$-language which is a ${\bf \Pi}^0_2$-complete subset 
of $X^\om$.
\nl
 One can prove that  $(g(W))^\om$  is 
${\bf \Si}^1_1$-complete, and hence a  non Borel set. This is done by reducing to this $\om$-language a well-known example of  a 
${\bf \Si}^1_1$-complete set: the set of infinite  binary trees labelled in the alphabet $\{0,1\}$  having  an infinite branch in the 
${\bf \Pi}^0_2$-complete set $(0^\star.1)^\om$. 
\nl More generally it is proved in  \cite[proof of Theorem 4.5 and Section 5]{Fin03a} 
that if $W^\om \subseteq X^\om$ is an $\om$-power which is ${\bf \Pi}^0_2$-hard, then the $\om$-power 
$(g(W))^\om \subseteq (X \cup \{d\})^\om$ is ${\bf \Si}^1_1$-complete,  and hence non Borel. 

\hs We use this result to prove our first proposition. In the sequel,  for two sets $A, B \subseteq 2^{X^\star}$  
we denote $A \leq B$ iff there is a continuous 
function $H: 2^{X^\star} \ra 2^{X^\star} $ such that $A = H^{-1}(B)$. So the relation $\leq $ is in fact the  Wadge reducibility  relation $\leq_W$.

\begin{prop}\label{prop}
The following relation holds :   $\mathcal{W}({\bf\Sigma}^0_{2}) \leq \mathcal{W}({\bf\Delta}_1^1).$
\end{prop}

\begin{proof}  We shall use the substitution $g$ defined above. Then let  $g' : X \cup \{d\} \ra \mathcal{P}(X^\star)$ be the substitution simply defined by 
$g'(0)=\{0\cdot 1\}$, $g'(1)=\{0\cdot 1^2\}$, and $g'(d)=\{0\cdot 1^3\}$. 
And let $G=g' \circ g$ be the substitution obtained by the composition of $g$ followed by $g'$. 
Then, for every dictionary $V \subseteq X^\star$, the language $G(V)$ is also a dictionary over the alphabet $X$ and $G(V^\om) = ( G(V) )^\om$. 
The substitution $G$ will provide the reduction $G: 2^{X^\star} \ra 2^{X^\star}$. 

\hs Firstly,   it is easy to see  that the mapping  $G: 2^{X^\star} \ra 2^{X^\star}$ is  continuous, \cite{Moschovakis80}.

\hs Secondly,  we claim that for every dictionary $V \subseteq X^\star$, it holds that:
   \[ V \in \mathcal{W}({\bf\Sigma}^0_{2}) 
\mbox{ if and only if } G(V) \in \mathcal{W}({\bf\Delta}_1^1). \]

\hs Assume first that  $ V \notin \mathcal{W}({\bf\Sigma}^0_{2})$. By definition of $\mathcal{W}({\bf\Sigma}^0_{2})$ 
this means that $V^\om$ is not a  ${\bf \Si}^0_2$-subset of $2^\om$. 
Then we can infer from Hurewicz's Theorem, see \cite[page 160]{Kechris94}, 
that the $\om$-power $V^\om$ is ${\bf \Pi}^0_2$-hard because it is an analytic subset of $2^\om$ 
which is not  a ${\bf \Si}^0_2$-set. 
Then it follows from \cite[proof of Theorem 4.5 and Section 5]{Fin03a} 
that the $\om$-power $(g(V))^\om \subseteq (X \cup \{d\})^\om$ is ${\bf \Si}^1_1$-complete, and  hence non Borel. 
It is now very easy to check, applying the second substitution $g'$, that the $\om$-power $(G(V))^\om \subseteq X^\om$ is also non Borel. 
This means that $G(V)$ does not belong to  the set $\mathcal{W}({\bf\Delta}_1^1)$. 

\hs Conversely assume now that $ V \in \mathcal{W}({\bf\Sigma}^0_{2})$. 
By definition of $\mathcal{W}({\bf\Sigma}^0_{2})$ this means that $V^\om$ is  a  ${\bf \Si}^0_2$-subset of $X^\om$, i.e.   
is  a countable union of closed sets $F_n \subseteq X^\om$, $n\geq 1$. 
Thus $V^\om = \bigcup_{n\geq 1} F_n$ and $G( V^\om ) = G( \bigcup_{n\geq 1} F_n ) = \bigcup_{n\geq 1}G( F_n ) $. 
\nl We are going to show that for every closed set $F \subseteq X^\om$, it holds that $G( F )$ is a Borel subset of $X^\om$. 

\hs Let then $F \subseteq X^\om$ be a closed set.  Then there is a tree $T \subseteq X^\star$ such that $F=[T]$, i.e. 
 $F$ is the set of the infinite branches of $T$. 
We first prove that $g(F)$ is Borel. For any $\om$-word $y \in (X \cup \{d\})^\om$, it holds 
that $y\in g(F)$ if and only if there exist $x\in F$ and sequences 
$u_i, v_i \in X^\star$, $i\geq 1$, such that :
\[y = x(1)\cdot (u_1\cdot d\cdot v_1)\cdot x(2)\cdot (u_2\cdot d\cdot v_2)\cdot x(3) \cdots \]
\noi where for each integer $i\geq 1$, ~$( |v_i|=2|u_i|)~~ or ~~( |v_i|=2|u_i|+1)$.
\nl Let then $T_1$ be the set of finite prefixes of such $\om$-words in the set $ g(F)$. The set $T_1 \subseteq (X \cup \{d\})^\star$ is a tree. 
We claim that $g(F)= [T_1] \cap (\{0, 1\}^\star\cdot d)^\om$. 

\hs The inclusion $g(F) \subseteq [T_1] \cap (\{0, 1\}^\star\cdot d)^\om$ is straightforward. 

\hs  To prove the inverse inclusion, let us consider an $\om$-word $x \in [T_1] \cap (\{0, 1\}^\star\cdot d)^\om$. 
\nl Then for each integer $n \geq 1$ there exists (at least) one finite sequence $(\varepsilon_i)_{1\leq i\leq n}\in \{0, 1\}^n$ 
and one finite word $a_{1}\cdot a_2  \cdots a_n \in X^\star$ 
and finite words $u_i$ and $v_i$ in $X^\star$, for $1 \leq i \leq n-1$,  and $u\in X^\star$, such that : 

\[a_1\cdot (u_1\cdot d\cdot v_1)\cdot a_2\cdot (u_2\cdot d\cdot v_2) \cdots a_{n-1}(u_{n-1}\cdot d\cdot v_{n-1})\cdot a_n\cdot u\cdot d  \sqsubseteq   x \]
\noi where for each integer $i \in [1, n]$, ~$( |v_i|=2|u_i|)$ iff $\varepsilon_i = 0$     and    $~( |v_i|=2|u_i|+1)$ iff $\varepsilon_i = 1$, 
and $a_1\cdot a_2 \cdots a_{n-1} \in T$. 

\hs Consider now all the ``suitable" sequences $(\varepsilon_i)_{1\leq i\leq n}\in \{0, 1\}^n$  defined as above. The set of these suitable sequences is closed 
under prefix. Therefore this set form a subtree of  $( \{0, 1\}^\star, \sqsubseteq)$,  which is finitely branching. 
This tree is infinite so by K\"onig's Lemma it has an infinite 
branch. therefore there exists an {\it infinite } sequence $(\varepsilon_i)_{1 \leq i < \om}\in \{0, 1\}^\om$ and 
one infinite word $a_{1}\cdot a_2\cdot  \cdots a_n \cdots  \in \Sio$ 
and finite words $u_i$ and $v_i$ in $X^\star$, for $1 \leq i <\om$,  such that : 

\[   x  = a_1\cdot (u_1\cdot d\cdot v_1)\cdot a_2\cdot (u_2\cdot d\cdot v_2) \cdots a_{n}(u_{n}\cdot d\cdot v_{n}) \cdots   \]
\noi where for each integer $i \geq 1$, ~$( |v_i|=2|u_i|)$ iff $\varepsilon_i = 0$     and    $~( |v_i|=2|u_i|+1)$ iff $\varepsilon_i = 1$, 
and $a_1\cdot a_2 \cdots a_{n} \cdots \in [T]=F$. 

\hs This shows that $x \in g(F)$. 

\hs Thus  $g(F)= [T_1] \cap (\{0, 1\}^\star\cdot d)^\om$ is the intersection of the closed set $[T_1] $ and 
of the ${\bf\Pi}^0_2$-set $(\{0, 1\}^\star\cdot d)^\om$. 
Then  $g(F)$ is a Borel ${\bf\Pi}^0_2$-set, and it is easy to see  that $G(F)$ is also Borel. 
\nl Assume now that  $V \in \mathcal{W}({\bf\Sigma}^0_{2})$, then $V^\om = \bigcup_{n\geq 1} F_n$, where  $F_n \subseteq X^\om$ are 
closed sets. Then   $G( V)^\om = G( V^\om ) = G( \bigcup_{n\geq 1} F_n ) = \bigcup_{n\geq 1}G( F_n ) $ is a Borel subset of $X^\om$, so 
$G(V)$  belongs to  the set  $\mathcal{W}({\bf\Delta}_1^1)$. 
\end{proof} 

\hs We can now improve the result : $\mathcal{W}({\bf\Delta}_1^1) \in {\bf\Sigma}^1_2(2^{2^{\star}})\!\setminus\! D_2(\boraone )$ proved in 
\cite{Lecomte-JSL}.

\begin{cor}
The following relation holds :  $\mathcal{W}({\bf\Delta}_1^1) \in {\bf\Sigma}^1_2(2^{2^{\star}})\!\setminus\! {\bf\Pi}^0_2$ 
\end{cor}

\begin{proof} It follows directly from the relations   $ \mathcal{W}({\bf\Sigma}^0_{2}) \in {\bf\Sigma}^1_2(2^{2^{\star}})\!\setminus {\bf\Pi}^0_2$ and 
$\mathcal{W}({\bf\Delta}_1^1) \in {\bf\Sigma}^1_2(2^{2^{\star}})$, proved by Lecomte 
in \cite{Lecomte-JSL},  and from Proposition \ref{prop} stating that $ \mathcal{W}({\bf\Sigma}^0_{2}) \leq \mathcal{W}({\bf\Delta}_1^1).$
\end{proof}

\begin{remark}{\rm 
We have obtained only a slight improvement of Lecomte's 
result that  $\mathcal{W}({\bf\Delta}_1^1) \in {\bf\Sigma}^1_2(2^{2^{\star}})\!\setminus\! D_2(\boraone )$. 
The question is  still open of the exact complexity of the two sets $\mathcal{W}({\bf\Delta}_1^1)$ and 
$\mathcal{W}({\bf\Sigma}_{2}^0)$ (and also of the other sets $\mathcal{W}({\bf\Sigma}^0_{k})$ and 
$\mathcal{W}({\bf\Pi}^0_{k})$). 
\nl However,  Proposition \ref{prop} could provide a better improvement of the lower bound on the complexity of $\mathcal{W}({\bf\Delta}_1^1)$ 
as soon as a better improvement of the lower bound on the complexity of $\mathcal{W}({\bf\Sigma}^0_{2})$ would be obtained. 
On the other hand, if one could  obtain a better upper bound on the complexity of the set $\mathcal{W}({\bf\Delta}_1^1)$, 
then this would provide, by Proposition  \ref{prop}, 
a better upper bound on the complexity of the set $\mathcal{W}({\bf\Sigma}_{2}^0)$. }
\end{remark}

\noi We consider now Borel $\om$-powers. It has been proved in \cite{Fin01a} that for each integer $n\geq 1$, there exist some 
$\om$-powers of (context-free) languages which are ${\bf\Pi}_{n}^0$-complete Borel sets. 
(We refer the reader for instance to \cite{ABB96} for definitions and properties of context-free languages). 
These results were obtained by the use of an operation 
$A \ra A^\approx$ over $\om$-languages  which is a  variant of Duparc's operation 
of exponentiation $A \ra A^\sim$, \cite{Duparc01}. 

\hs We first recall  the definition of the operation  $A \ra A^\sim$. Notice that this operation is defined over sets of finite {\it or } infinite words, called 
{\it conciliating sets} in \cite{Duparc01}. 

\begin{defn}[Duparc \cite{Duparc01}]\label{til}
Let  $X$ be a finite alphabet,  $\tla  \notin X$, and  
 $x$ be a finite or infinite word over the alphabet $Y=X\cup \{\tla\}$.
\nl Then  $x^\tla$ is inductively defined by:
\nl $\lambda^\tla =\lambda$,
\nl and for a finite word $u\in (X \cup \{\tla\})^\star$:
\nl $(u\cdot a)^\tla=u^\tla\cdot a$, if $a\in X$,
\nl $(u\cdot \tla)^\tla =u^\tla$  with its last letter removed if $|u^\tla|>0$,
\nl i.e. $(u\cdot \tla)^\tla =u^\tla(1)\cdot u^\tla(2)\cdots u^\tla(|u^\tla|-1)$  if
$|u^\tla|>0$,
\nl $(u\cdot \tla)^\tla=\lambda$  if $|u^\tla|=0$,
\nl and for $u$ infinite:
\nl $(u)^\tla = \lim_{n\in\om} (u[n])^\tla$, where, given $\beta_n$ and $v$
in   $X^\star$,
\nl $v\sqsubseteq \lim_{n\in\om} \beta_n \lra  \exists n \fa p\geq n\quad
\beta_p[|v|]=v$.
\nl(The finite {\it or} infinite word $\lim_{n\in\om} \beta_n$ is
determined by the set of its (finite) prefixes).
\end{defn}

\begin{remark}{\rm 
For $x \in Y^{\leq \om}$, $x^\tla$ denotes the string $x$, once every $\tla$
occuring in $x$
has been ``evaluated" to the back space operation, 
proceeding from left to right inside $x$. In other words $x^\tla = x$ from
which every
 interval of the form $``a\tla "$ ($a\in X$) is removed. The letter  $\tla$ may be called an ``eraser".} 
\end{remark}

\noi For example if $u=(a\tla)^n$, for $n$ an integer $\geq 1$, or
$u=(a\tla)^\om$,  or $u=(a\tla\tla)^\om$,  then $(u)^\tla=\lambda$.
If $u=(ab\tla)^\om$ then $(u)^\tla=a^\om$  and
 if $u=bb(\tla a)^\om$ then $(u)^\tla=b$.

\hs Let us notice that in Definition \ref{til} the limit is not defined in
the usual way:
\nl for example if $u=bb(\tla a)^\om$ the finite word  $u[n]^\tla$ is
alternatively
equal to $b$ or to $ba$: more precisely $u[2n+1]^\tla=b$ and
$u[2n+2]^\tla=ba$ for every
integer $n\geq 1$ (it holds also that
$u[1]^\tla=b$ and  $u[2]^\tla=bb$). Thus Definition \ref{til} implies that
$\lim_{n\in\om} (u[n])^\tla = b$ so $u^\tla=b$.

\hs  We can now define the operation $A \ra A^\sim$ of
{\it exponentiation of conciliating sets}:

\begin{defn}[Duparc \cite{Duparc01}]
For $A\subseteq X^{\leq \om}$ and $\tla \notin X$, let 
\[A^\sim =_{df} \{x\in (X\cup \{\tla\})^{\leq \om} \mid  x^\tla\in
A\}.\]
\end{defn}

\noi We now define the variant $A \ra A^\approx$ of the  operation 
 $A \ra A^\sim$. 

\begin{defn}[\cite{Fin01a}]\label{til2}
Let  $X$ be a finite alphabet,  $\tla  \notin X$, and 
 $x$ be a finite or infinite word over the alphabet $Y=X\cup \{\tla\}$.
\nl Then  $x^{\hookleftarrow}$ is inductively defined by:
\nl $\lambda^{\hookleftarrow} =\lambda$,
\nl and for a finite word $u\in (X \cup \{\tla\})^\star$:
\nl $(u\cdot a)^{\hookleftarrow}=u^{\hookleftarrow}\cdot a$, if $a\in X$,
\nl $(u\cdot \tla)^{\hookleftarrow} =u^{\hookleftarrow}$  with its last letter removed if $|u^{\hookleftarrow}|>0$,
\nl $(u\cdot \tla)^{\hookleftarrow}$ is undefined if $|u^{\hookleftarrow}|=0$,
\nl and for $u$ infinite:
\nl $(u)^{\hookleftarrow} = \lim_{n\in\om} (u[n])^{\hookleftarrow}$, where, given $\beta_n$ and $v$
in   $X^\star$,
\nl $v\sqsubseteq \lim_{n\in\om} \beta_n \lra  \exists n \fa p\geq n\quad
\beta_p[|v|]=v$.
\end{defn}

\noi The difference between the definitions  of $x^\tla$ and $x^{\hookleftarrow}$ is that here we have 
 added the convention that $(u.\tla)^{\hookleftarrow}$ is undefined if $|u^{\hookleftarrow}|=0$, i.e. 
when the last letter $\tla$ can not be used as an eraser (because every letter 
of $X$ in $u$ has already been erased by some erasers $\tla$ placed in $u$).  
 For example if $u=\tla (a\tla)^\om$ or $u=a\tla\tla a^\om$ or $u=(a \tla\tla)^\om$,  then $(u)^{\hookleftarrow}$ is undefined. 

\begin{defn}
For $A\subseteq X^{\leq \om}$, ~~  $A^\approx = \{x\in (X\cup \{\tla\})^{\leq \om} \mid  x^{\hookleftarrow}       \in
A\}$.
\end{defn}

\noi The operation $A \ra A^\sim$ was used by Duparc in his study of the Wadge hierarchy, \cite{Duparc01}.  The result 
stated in the following lemma will be  important in the sequel. 

\begin{lem}\label{sim-approx}
Let $X$ be a finite alphabet and $L \subseteq X^\om$. Then the two $\om$-languages $L^\sim$ and $L^\approx$ are Wadge equivalent, i.e. 
$L^\sim \equiv_W L^\approx$. 
\end{lem}

\begin{proof}  Let $X$ be a finite alphabet and $L \subseteq X^\om$. 
We are going to prove that $L^\sim \equiv_W L^\approx$, using Wadge games. 

\begin{enumerate}
\ite[a)] In the  Wadge game $W(L^\sim, L^\approx)$ 
  the player in charge of 
$L^\approx$ has clearly a winning strategy 
which consists in copying the play of the other player except if player $1$ writes the eraser
$\tla$ but he has nothing to erase. In this case player $2$ writes for example a letter 
$a\in X$ and the eraser $\tla$ at the next step of the play.
Now  if, in $\om$ steps, player $1$ has written the $\om$-word $\alpha$ and 
player $2$ has written the $\om$-word $\beta$, then it is easy to see that 
$[\alpha^\tla = \beta^{\hookleftarrow}]$ and then $\alpha \in  L^\sim$ iff $\beta \in L^\approx$. 
Thus player $2$ has 
a winning strategy in the Wadge game $W(L^\sim, L^\approx)$.
 
\ite[b)] Consider now the Wadge game $W(L^\approx , L^\sim)$. 
The only extra  
possibility for  player $1$  in charge of $L^\approx$ is to get out of the set $L^\approx$ by 
writing the eraser $\tla$ when in fact there is no letter of his previous play to erase.
But then his final play is surely outside $L^\approx$.  
If this happens at some point of the play, then player $2$ may writes the eraser $\tla$ forever. Then,  after $\om$ steps, player $2$ has written an 
infinite word $\beta$ such that $\beta^\tla = \lambda$. In particular,  $\beta^\tla \notin L$ because $\beta^\tla$ is not an infinite word, and $\beta \notin L^\sim$. 
On the other hand player $1$ has written an infinite word $\alpha$ such that $\alpha^{\hookleftarrow}$ is undefined, hence $\alpha \notin L^\approx$. 
Therefore player $2$ wins the play in this case too, and 
player 2 has a winning strategy in the Wadge game $W(L^\approx , L^\sim)$.   

\end{enumerate}

\end{proof}

\noi The operation $A \ra A^\approx$ is very useful  in the study of $\om$-powers because it can be defined with the notion of substitution and 
preserves the $\om$-powers of finitary languages. 
   Let $L_1 = \{ w\in  (X \cup\{\tla\})^\star \mid w^{\hookleftarrow}=\lambda \}$.  $L_1$ is a 
context free (finitary) 
language generated by the context free grammar with the following productions:
~~  $(S, ~ aS\tla S)$  with $a\in X$; 
and  $(S, ~ \lambda)$. 
\nl 
Then, for each $\om$-language $A \subseteq X^{\om}$, the  $\om$-language 
$A^\approx \subseteq (X \cup\{\tla\})^\om$
is obtained by substituting in $A$ the language $L_1\cdot a$ for each letter $a\in X$. 
This implies that  the  operation $A \ra A^\approx$ preserves the $\om$-powers 
of finitary languages. This is stated in the following lemma. 

\begin{lem}[\cite{Fin01a}]
Let $X$ be a finite alphabet and let $h$ be the 
substitution defined by 
$h(a)=L_1\cdot a$ for every letter $a\in X$. 
\nl If 
$A=V^\om$ for some language $V   \subseteq X^\star$, then 
$A^\approx=h(V^\om)=( h(V) )^\om$. Thus, if $A$ is an $\om$-power, then  $A^\approx$   is also an $\om$-power. 
\end{lem}

\hs We now recall  the operation $A \ra A^b$  used by Duparc in his study of the Wadge hierarchy, \cite{Duparc01}. 
  For $A\subseteq X^{\leq \om}$ and $b$ a letter not in $X$, $A^b$ is the $\om$-language over  $X \cup \{b\}$  
which is defined by : 
 \[A^b = \{ x\in (X \cup \{b\})^\om \mid x( /b)\in A \}\]
\noi where $x( /b)$ is the sequence obtained from $x$ 
when removing every occurrence of the letter $b$.

\hs We can now state the following lemma. 

\begin{lem}\label{b}  Let $X$ be a finite alphabet having at least two elements  and $A \subseteq X^\om$. 
\begin{enumerate}
\ite For each integer $k\geq 2$, ~~~~  $A$ is a ${\bf\Pi}^0_k$-subset of $X^\om$ iff $A^b$ is a ${\bf\Pi}^0_k$-subset of $(X \cup \{b\})^\om$. 
\ite For each integer $k\geq 3$, ~~~~  $A$ is a ${\bf\Si}^0_k$-subset of $X^\om$ iff $A^b$ is a ${\bf\Si}^0_k$-subset of $(X \cup \{b\})^\om$. 
\end{enumerate}
\end{lem}

\begin{proof}  We denote by $Z^\infty$ the set of infinite words 
in $(X \cup \{b\})^\om$ having infinitely many letters in $X$. The set $Z^\infty=\{ x \in (X \cup \{b\})^\om \mid x( /b)\in X^\om \}$
 is a well known example of ${\bf\Pi}^0_2$-subset of $(X \cup \{b\})^\om$, 
\cite{Kechris94,PerrinPin}.  Notice that $Z^\infty$, equipped with the induced topology, 
 is a topological subspace of the Cantor space $(X \cup \{b\})^\om$. 
One can define the Borel hierarchy on the topological space 
$Z^\infty$ as in the case of the Cantor space, see \cite[page 68]{Kechris94}.  Then one can prove by induction that,  
for each non-null  countable ordinal  $\alpha$,  the 
${\bf \Si}^0_\alpha$ (respectively,  ${\bf \Pi}^0_\alpha$)-subsets   of   $Z^\infty$  are the intersections of   
${\bf \Si}^0_\alpha$ (respectively,  ${\bf \Pi}^0_\alpha$)-subsets of  $(X \cup \{b\})^\om$ with the set  $Z^\infty$, see 
 \cite[page 167]{Kechris94}. 

\hs 
Let now $\phi$ be the function from 
$Z^\infty$ into $X^\om$ defined by $\phi(x) = x( /b)$. It is easy to see  that, for each $A \subseteq X^\om$, it holds that 
$\phi^{-1}(A)=A^b$. On the other hand, 
 the function $\phi$ is continuous. Thus the inverse image of an open (respectively, closed) subset of $X^\om$ is an 
open (respectively, closed) subset of  $Z^\infty$. And one can  prove by induction that, for each non-null   countable ordinal  $\alpha$,  
the inverse image of a ${\bf \Si}^0_\alpha$ (respectively,  ${\bf \Pi}^0_\alpha$)-subset   of  $X^\om$ is
 a ${\bf \Si}^0_\alpha$ (respectively,  ${\bf \Pi}^0_\alpha$)-subset   of   $Z^\infty$, i.e. the intersection of  a  
${\bf \Si}^0_\alpha$ (respectively,  ${\bf \Pi}^0_\alpha$)-subset   of   $(X \cup \{b\})^\om$ with the set  $Z^\infty$. 

\hs 
 Let now $k\geq 2$ and $A \subseteq X^\om$ be a ${\bf\Pi}^0_k$-subset of $X^\om$.  Then $\phi^{-1}(A)=A^b$ is a ${\bf\Pi}^0_k$-subset of 
of   $Z^\infty$, i.e. the  intersection of  a   ${\bf \Pi}^0_k$-subset   of   $(X \cup \{b\})^\om$ with the set  $Z^\infty$. 
But $Z^\infty$ is a ${\bf\Pi}^0_2$-subset of $(X \cup \{b\})^\om$ thus $\phi^{-1}(A)=A^b$ is the intersection of two ${\bf \Pi}^0_k$-subsets of 
 $(X \cup \{b\})^\om$, hence also a ${\bf \Pi}^0_k$-subset of  $(X \cup \{b\})^\om$. 

\hs 
 In a similar way we prove that if $k\geq 3$ and $A \subseteq X^\om$ is a ${\bf\Si}^0_k$-subset of $X^\om$, then $\phi^{-1}(A)=A^b$ is a 
${\bf\Si}^0_k$-subset of  $(X \cup \{b\})^\om$. 

\hs  Conversely assume that for some integer $k\geq 2$ and $A \subseteq X^\om$ the set $A^b$ is a ${\bf \Pi}^0_k$-subset of  $(X \cup \{b\})^\om$. 
Notice that $X^\om$ is a closed subset of $(X \cup \{b\})^\om$. 
Thus $A=A^b \cap X^\om$ is the intersection of two ${\bf \Pi}^0_k$-subsets of  $(X \cup \{b\})^\om$, hence also a 
${\bf \Pi}^0_k$-subset of $(X \cup \{b\})^\om$. And $A=A\cap X^\om$ so $A$ is also a ${\bf \Pi}^0_k$-subset of $X^\om$. 

\hs  In a similar way we prove that if for some integer $k\geq 3$ and $A \subseteq X^\om$ the set $A^b$ is a ${\bf \Si}^0_k$-subset of  $(X \cup \{b\})^\om$, 
then $A$ is also a ${\bf \Si}^0_k$-subset of  $X^\om$. 
\end{proof} 

\begin{lem}\label{approx}
 Let $X$ be a finite alphabet having at least two elements  and $A \subseteq X^\om$. 
\begin{enumerate}
\ite For each integer $k\geq 3$, ~~  $A$ is a ${\bf\Si}^0_k$-subset of $X^\om$ iff $A^\approx$ is a ${\bf\Si}^0_{k+1}$-subset of $(X \cup \{\tla\})^\om$. 
\ite For each integer $k\geq 2$, ~~  $A$ is a ${\bf\Pi}^0_k$-subset of $X^\om$ iff $A^\approx$ is a ${\bf\Pi}^0_{k+1}$-subset of $(X \cup \{\tla\})^\om$. 
\end{enumerate}
\end{lem}

\begin{proof}  Let $X$ be a finite alphabet having at least two elements,    $A \subseteq X^\om$,  and  $k\geq 3$ be an integer.
 Then  the following equivalences hold: 

\hs ~~~~~~~~$A \in $ ${\bf\Si}^0_k$
\nl $\longleftrightarrow$  ~$A^b \in $ ${\bf\Si}^0_k$   ~~~~~~~~~~~~~~~~~~ by Lemma \ref{b}. 
\nl $\longleftrightarrow$  ~$A^b \leq_W B^b$ ~~~~~~~~~~for some $B\subseteq X^{\leq \om}$ such that $B^b$ is ${\bf\Si}^0_k$-complete. 
\nl $\longleftrightarrow$  ~$(A^\sim)^b \leq_W (B^\sim)^b$  ~~~~~~~~~~~~ by \cite[Proposition 23]{Duparc01}.
\nl $\longleftrightarrow$  ~$(A^\sim)^b \in $ ${\bf\Si}^0_{k+1}$,  because $(B^\sim)^b$ is ${\bf\Si}^0_{k+1}$-complete by \cite[Lemma 31]{Duparc01}.
\nl $\longleftrightarrow$  ~$A^\sim \in $ ${\bf\Si}^0_{k+1}$  ~~~~~~~~~~~~ ~by Lemma \ref{b}.
\nl $\longleftrightarrow$  ~$A^\approx \in $ ${\bf\Si}^0_{k+1}$ ~~~~~~~~~~~~~~by Lemma \ref{sim-approx}.

\hs In a very similar way we prove that if $k\geq 2$ is an integer, then $A \in $ ${\bf\Pi}^0_k$ iff $A^\approx \in $ ${\bf\Pi}^0_{k+1}$. 
\end{proof} 

\hs We now state the following result about the classes $\mathcal{W}({\bf\Pi}^0_{k})$. 

\begin{prop}
For each   integer $k \geq 2$ it  holds that:  $\mathcal{W}({\bf\Pi}^0_k) \leq \mathcal{W}({\bf\Pi}^0_{k+1}).$
\end{prop}

\begin{proof}  We shall use the substitution $h$ defined above. 
\nl Let then $h' : \{0, 1, \tla\} \ra \mathcal{P}(\{0, 1\}^\star)$ be the substitution simply defined by 
$h'(0)=\{0\cdot 1\}$, $h'(1)=\{0\cdot 1^2\}$, and $h'(\tla)=\{0\cdot 1^3\}$. 
And let $H=h' \circ h$ be the substitution obtained by the composition of $h$ followed by $h'$. 
Then, for every dictionary $V \subseteq X^\star =\{0, 1\}^\star$, 
the language $H(V)$ is also a dictionary over the alphabet $X$ and $H(V^\om) = ( H(V) )^\om$. 
The substitution $H$ will provide the reduction $H: 2^{X^\star} \ra 2^{X^\star}$. 

\hs It is easy to see  that the mapping  $H: 2^{X^\star} \ra 2^{X^\star}$ is continuous, \cite{Moschovakis80}.

\hs We claim that for every dictionary $V \subseteq X^\star$, it holds that   $ V \in \mathcal{W}({\bf\Pi}^0_k)$ if and only if 
$H(V) \in \mathcal{W}({\bf\Pi}^0_{k+1})$. 

\hs Firstly by definition of the class $\mathcal{W}({\bf\Pi}^0_k)$ it holds that for every dictionary $V \subseteq X^\star$, $V$ is in the class 
$\mathcal{W}({\bf\Pi}^0_k)$ iff 
 $V^\om$ is a  ${\bf\Pi}^0_k$-set. 
By Lemma \ref{approx},   $V^\om$ is a  ${\bf\Pi}^0_k$-set 
 iff $(V^\om)^\approx$ is a ${\bf\Pi}^0_{k+1}$-set. But $(V^\om)^\approx = h(V^\om)=h(V)^\om$. 
Thus $V$ is in the class $\mathcal{W}({\bf\Pi}^0_k)$ iff 
$h(V)^\om$  is in the class ${\bf\Pi}^0_{k+1}$. It is now easy to see, using the coding $h'$ that this is equivalent to the assertion ``$( H(V) )^\om$ is in 
the class ${\bf\Pi}^0_{k+1}$", i.e. $H(V) $ is in the class $\mathcal{W}({\bf\Pi}^0_{k+1})$. 
\end{proof} 

\hs In a very similar manner, we can prove the following result about the sets $\mathcal{W}({\bf\Si}^0_k)$ for integers $k\geq 3$.

\begin{prop}
For each   integer $k \geq 3$ it  holds that:   $\mathcal{W}({\bf\Si}^0_k) \leq \mathcal{W}({\bf\Si}^0_{k+1}).$
\end{prop}

\begin{remark}{\rm 
Notice that here $k\geq 3$ because 
for $L \subseteq X^\om$ then $L$ may be  in the class ${\bf\Si}^0_2$ while  
$L^b \subseteq (X \cup \{b\})^\om$ is not in the class ${\bf\Si}^0_2$.  For instance $L=\{0, 1\}^\om \subseteq \{0, 1\}^\om$  
is open and closed hence also  in the class ${\bf\Si}^0_2$. But the $\om$-language $L^b$ is simply the set of $\om$-words over the alphabet $\{0, 1, b\}$ 
which contain infinitely many letters $0$ or $1$ and it is a ${\bf\Pi}^0_2$-complete, hence non ${\bf\Si}^0_2$,  subset of $\{0, 1, b\}^\om$.  }
\end{remark}

\section{ Concluding remarks} 

\noi 
Lecomte proved  that    for every countable ordinal  $\xi \geq 2$   (respectively, $\xi \geq 3$),  
$\mathcal{W}({\bf\Pi}^0_\xi) \in {\bf\Sigma}^1_2(2^{2^{\star}})\!\setminus\! D_2(\boraone )$  (respectively, 
$\mathcal{W}({\bf\Si}^0_\xi) \in {\bf\Sigma}^1_2(2^{2^{\star}})\!\setminus\! D_2(\boraone )$). 
Finkel and Lecomte proved that for every countable ordinal  $\xi \geq 3$,
the sets $\mathcal{W}({\bf\Pi}^0_\xi)$ and   $\mathcal{W}({\bf\Si}^0_\xi)$ are actually ${\bf\Pi}^1_1$-hard. 
The exact complexity of the  sets 
$\mathcal{W}({\bf\Pi}_{\xi}^0)$ and   $\mathcal{W}({\bf\Si}_\xi^0)$ is still unknown, 
 but our new results could help to determine it. 

\hs {\bf  Acknowledgements.} Thanks to  Dominique Lecomte and  to the anonymous referee for useful comments on a preliminary version 
of this paper.

\providecommand{\WileyBibTextsc}{}
\let\textsc\WileyBibTextsc
\providecommand{\othercit}{}
\providecommand{\jr}[1]{#1}
\providecommand{\etal}{~et~al.}


\begin{thebibliography}{[10]}

\bibitem{Lecomte-JSL}
 \textsc{D.~Lecomte},
 \jr{Journal of Symbolic Logic} \textbf{70}(4), 1210--1232 (2005).


\othercit
\bibitem{Staiger97}
 \textsc{L.~Staiger},
$\omega$-languages,
 in: Handbook of formal languages, Vol.\ 3,  (Springer, Berlin, 1997),
  pp.\,339--387.


\othercit
\bibitem{Niwinski90}
 \textsc{D.~Niwinski},
A problem on $\om$-powers,
 in: 1990 Workshop on Logics and Recognizable Sets,  (University of Kiel,
  1990).


\othercit
\bibitem{Simonnet92}
 \textsc{P.~Simonnet},
Automates et th\'eorie descriptive,
PhD thesis, Universit\'e Paris VII, 1992.


\othercit
\bibitem{Staiger97b}
 \textsc{L.~Staiger},
On $\omega$-power languages,
 in: New Trends in Formal Languages, Control, Coperation, and Combinatorics, 
  Lecture Notes in Computer Science,  Vol.\,1218 (Springer-Verlag, Berlin,
  1997),  pp.\,377--393.


\bibitem{Fin03a}
 \textsc{O.~Finkel},
 \jr{Theoretical Computer Science} \textbf{290}(3), 1385--1405 (2003).


\bibitem{Fin01a}
 \textsc{O.~Finkel},
 \jr{Theoretical Computer Science} \textbf{262}(1--2), 669--697 (2001).


\bibitem{Fin04-FI}
 \textsc{O.~Finkel},
 \jr{Fundamenta Informaticae} \textbf{62}(3--4), 333--342 (2004).


\othercit
\bibitem{Fin-Dup}
 \textsc{J.~Duparc} and  \textsc{O.~Finkel},
An $\omega$-power of a context free language which is {B}orel above
  ${\Delta}_\omega^0$,
 in: Proceedings of the International Conference Foundations of the Formal
  Sciences V : Infinite Games, November 26th to 29th, 2004, Bonn, Germany, 
  Studies in Logic, College Publications at King's College,  Vol.\,11 (London,
  2007),  pp.\,109--122.


\othercit
\bibitem{Fin-Lec}
 \textsc{O.~Finkel} and  \textsc{D.~Lecomte},
There exist some $\omega$-powers of any {B}orel rank,
 in: Proceedings of the 16th EACSL Annual International Conference on Computer
  Science and Logic, CSL 2007, Lausanne, Switzerland, September 11-15, 2007, 
  Lecture Notes in Computer Science,  Vol.\,4646 (Springer, 2007),
  pp.\,115--129.


\bibitem{Fink-Lec2}
 \textsc{O.~Finkel} and  \textsc{D.~Lecomte},
 \jr{Annals of Pure and Applied Logic} \textbf{160}(2), 163--191 (2009).


\othercit
\bibitem{Thomas90}
 \textsc{W.~Thomas},
Automata on infinite objects,
 in: Handbook of Theoretical Computer Science, edited by J.~van Leeuwen
  (Elsevier, 1990),  pp.\,135--191.


\othercit
\bibitem{Moschovakis80}
 \textsc{Y.\,N. Moschovakis},
Descriptive set theory (North-Holland Publishing Co., Amsterdam, 1980).


\othercit
\bibitem{LescowThomas}
 \textsc{H.~Lescow} and  \textsc{W.~Thomas},
Logical specifications of infinite computations,
 in: A Decade of Concurrency, edited by J.\,W. de~Bakker, W.\,P. de~Roever,
  and G.~Rozenberg, Lecture Notes in Computer Science Vol.\,803 (Springer,
  1994),  pp.\,583--621.


\othercit
\bibitem{Kechris94}
 \textsc{A.\,S. Kechris},
Classical descriptive set theory (Springer-Verlag, New York, 1995).


\othercit
\bibitem{PerrinPin}
 \textsc{D.~Perrin} and  \textsc{J.\,E. Pin},
Infinite words, automata, semigroups, logic and games, Pure and Applied
  Mathematics,  Vol.\,141 (Elsevier, 2004).


\othercit
\bibitem{Wadge83}
 \textsc{W.~Wadge},
Reducibility and determinateness in the Baire space,
PhD thesis, University of California, Berkeley, 1983.


\othercit
\bibitem{ABB96}
 \textsc{J.\,M. Autebert},  \textsc{J.~Berstel},  and  \textsc{L.~Boasson},
Context free languages and pushdown automata,
 in: Handbook of formal languages, Vol.\ 1,  (Springer-Verlag, 1996).


\bibitem{Duparc01}
 \textsc{J.~Duparc},
 \jr{Journal of Symbolic Logic} \textbf{66}(1), 56--86 (2001).


\end{thebibliography}
\end{document}